\documentstyle[12pt,amscd]{amsart}
\newcommand{\C}{\Bbb{C}}
\newcommand{\pts}{\hbox{\rm{pts}}}
\newcommand{\pr}{\bold P}
\newcommand{\LL}{\cal{P}}
\newcommand{\VV}{\cal{V}}
\newcommand{\SS}{\cal{S}}

\newtheorem{thm}{Theorem}[section]
\newtheorem{cor}[thm]{Corollary}
\newtheorem{lem}[thm]{Lemma}
\newtheorem{pro}[thm]{Proposition}

\newtheorem{defn}{Definition}

\newtheorem{rem}[thm]{Remark}

\begin{document}
\title {SEVERI DEGREES IN COGENUS 4}
\author {Youngook Choi}
\date {January 16, 1996}
\maketitle

\pagestyle{plain}
\section{Summary}
Let $\LL (d)$
 be the linear system of degree $d$ curves in $\pr^2$.
Then $\LL (d)$ is a projective space of dimension $\binom{d+2}{2}-1$.
The Severi variety $\VV (\delta, d) \subset \LL (d)$ is the subset
corresponding to reduced, nodal curves with $\delta$ nodes.
The ``well-known'' problem, whether $\VV (\delta, d)$ is irreducible
or not, is solved affirmatively by Harris [H], and Ran [R1]. The next question
about the  Severi variety is to find its degree, $\text{N}(\delta,d)$.
 In the paper [R2],
 Ran gave the recursive formulae. In this paper, we will give closed-form
formulae for cogenus 3 and 4 cases using his method.
 These formulae coincide with
those of I. Vainsencher [V], and
 for cogenus 3 case, that of J. Harris and R. Pandharipande [H-P].
Ran gave the formula of cogenus 2 case using his method in the
paper [R3], namely,
$$
\text{N}(2,d) = 3/2(d-1)(d-2)(3d^2-3d-11)
$$
 Therefore, this paper is an extension of that paper.
 One small
benefit of our approach is that we can calculate, as a special case, the
formulae of degree of the locus of all curves having less
than three nodes and some
tangency conditions to fixed line.

Another result of this paper is that we calculate the degree of
the polynomial N($\pi,\delta, d$) in $d$, which is
the degree of the locus of curves $C$ in $\pr^2$ with degree $d$
and having $\delta$ nodes and $C \cap \text{L}$ is
{\it of type $\pi$} to fixed line L. Using this result, we also
calculate the coefficients of two leading terms
of Severi polynomial, N($\delta, d$). The results are :
for any $\delta$ and any partition $\pi = [\ell_1,...,\ell_n]$,
$$
\begin{array}{ll}
\text{deg} \ \text{N} (\pi,\delta,d) = 2\delta + \sum_{i=2}^n \ell_i
& (\text{Prop. 4.2}) \\
a_{2\delta}^{\delta} = {3^{\delta}\over{\delta !}}
& (\text{Cor. 4.3}) \\
a_{2\delta - 1}^{\delta} = -2\delta a_{2\delta}^{\delta}
                         ={ {-2 \times 3^{\delta}}\over{(\delta-1)!}}
& (\text{Cor. 4.5}) \ ,
\end{array}
$$
\noindent
 where $a_{\beta}^{\delta}$ is the coefficient of the degree $\beta$ term
of the polynomial N($\delta,d$). So,
$$
\text{N} (\delta,d) =
{3^{\delta}\over{\delta !}}d^{2\delta} -{{ 2 \times 3^{\delta}}
\over{(\delta -1)!}}{d^{2\delta - 1}} + \text{lower degree terms.}
$$
For $\delta \leq $ 6, these two coefficients coincide with those of the
polynomials of I. Vainsencher [V].

This paper is written under the guidence of Z. Ran. Without his advice, it
would be
impossible for me to write  this paper. I deeply appreciate his advice.

\section{Set of divisors of $\pr^1$}
In this section, we calculate the degree of set
of divisors of $\pr^1$ because we need the degree of the locus of smooth curves
with tangency conditions to fixed line.
\begin{defn} A partition $\pi = [\ell_1,...,\ell_n]$ is a sequence
of nonnegative integers.
\end{defn}
\begin{defn} A divisor $D$ of $\pr^1$ is of type $\pi$
if it has the form $D=\sum_{i=1}^n \sum_{j=1}^{\ell_i}iP_{ij}$
for some distinct points $P_{ij}$.
\end{defn}
\begin{rem}

\begin{enumerate}
\item Degree of a divisor $D$ of type $\pi$ is
$\sum_{i=1}^n i\ell_i$.
\item Let $\Gamma_\pi$ be the closure of the set of all divisors
in $\pr^1$ of type $\pi$. If we identify the space of all divisors of degree
$N$ on $\pr^1$ with $\pr^N$, then $\Gamma_\pi$ is a subvariety of $\pr^N$,
where $N=\sum_{i=1}^n i\ell_i$.
\item $\text{\rm dim} \ \Gamma_\pi = \sum_{i=1}^n \ell_i$.
\item $m(\pi) \buildrel \rm def \over = \prod i^{\ell_i}$
\item $n(\pi) \buildrel \rm def \over =
{(\sum \ell_i)! \over \ell_1! \dotsm \ell_n!}$
\end{enumerate}
\end{rem}
\begin{lem} $\text{\rm deg} \ \Gamma_\pi=m(\pi)n(\pi)$
\end{lem}
{\it Proof} \ \ \
Let $Q \in \pr^1$, and $\text{N} = \sum_{i=1}^{n} i\ell_i $.
Let $H_Q=\{ D| D \ni Q, \ \text{\rm deg} \ D = \text{N}\}$. Then
$H_Q$ is a hyperplane in $\pr^N$. Then
$\Gamma_\pi \cap H_Q = \sum_{i=1}^n i\Gamma_{\pi, i}$, where
$\Gamma_{\pi, i}= \{ D \in \Gamma_{\pi}| D \ni iQ, \ D \not\ni (i+1)Q \}$.
 Therefore,
$\Gamma_{\pi} \cap H_{Q_1} \cap \dotsb \cap H_{Q_{\text{L}} }$,
where $\text{L}=\sum \ell_i$,
has $n(\pi)$ points (set theoretically) and each point has multiplicity
$m(\pi)$. Therefore, $deg \ \Gamma_{\pi}=n(\pi)m(\pi)$.
\begin{cor} Assume that $L$ is a line in $\pr^2$.
Fix a partition $\pi=[a_1, \dotsc , a_n].$ Let $\SS (\pi)$ be
the (locally closed) set of all smooth curves $C$ of degree $d$ in $\pr^2$
such that $C \cap \text{L}$ is of type $\pi$,
Then $deg \ \overline {\SS (\pi)} = n(\pi)m(\pi)$.
\end{cor}

\section{Recursion method}
In this section, we want to describe the method of calculating the
degree of $\overline{\VV (\delta, d)}$. You can also consult
 the paper [R2] or [R3].

\subsection{Degeneration of $\pr^2$}

Let $S$ be the blow up of $\C \times \pr^2$ at a point $(0, p)$,
and let $\pi^\prime : S \rightarrow \C \times \pr^2$ be the
blowdown map and $\pi = \rho_2 \circ \pi^\prime : S \rightarrow \C$
be the composition map. Then
\begin{equation} \pi^{-1} (t) =
 \begin{cases}
  \pr^2& \text{if $t \not= 0$}\\
  \pr^2 \cup \tilde{\pr^2}& \text{if $t=0$} \ ,
 \end{cases}
\end{equation}
\noindent
where $\tilde{\pr^2}$ is the blow up of $\pr^2$ at $p$.
In the case $\pi^{-1} (0) = \pr^2 \cup \tilde{\pr^2}, \
\text{E} \buildrel \rm def \over =
\pr^2 \cap
\tilde{\pr^2} $ called the axis is a line in $\pr^2$,
 and an exceptional divisor
in $\tilde{\pr^2}$. This gives the degeneration of $\pr^2$ to
$\pr^2 \cup \tilde{\pr^2}$.

\subsection{Degeneration of curve in $\pr^2$}

As $\pr^2$ is degenerated to $\pr^2 \cup \tilde{\pr^2}$,
 a curve in $\pr^2$ also degenerates to
$C_1 \cup C_2$, where $C_1 \subset \pr^2$ and $C_2 \subset \tilde{\pr^2}$.
The point is that the degenerated curve $C_1$
has smaller degree, therefore we can apply
recursion and the other part
$C_2$ has simple configuration, so it's easy to
count its dimension, degree etc.
 Clearly, $\text{N} (\delta,d)
\buildrel \rm def \over = \deg \ \overline{\VV (\delta, d)} $
  is the same as the number of  plane curves
with $\delta$ nodes and
 through
 $(\binom{d+2}{2}-1-\delta)$ generally given points.
Let's put
 $\binom{d+1}{2}-1-\delta$ points into
$\pr^2$ and  $(d+1)$ points into $\tilde{\pr^2}$, i.e,
$$
\begin{array}{ccccc}
\pr^2 & \sim\longrightarrow & \pr^2 & \cup & \tilde{\pr^2} \\
\cup  &                 & \cup  &      & \cup          \\
C     & \sim\longrightarrow & C_1   & \cup & C_2          \\
\cup  &                 & \cup  &      & \cup          \\
\binom{d+2}{2}-1-\delta \ \pts
 & \sim\longrightarrow & \binom{d+1}{2}-1-\delta \ \pts &  & (d+1) \ \pts \ ,
\end{array}
$$
\noindent
 and think of all possible kinds
of degeneration of $\delta$ nodes.
Then, $\text{\rm deg} \ C_1 = d-1$, and $C_2$ is a smooth rational
curve plus several rulings (possibly multiple) [R3].
After all, the number of all possible curves
 in $\pr^2 \cup \tilde{\pr^2}$
satisfying given conditions is $ \text{N} (\delta,d)$.

\subsection{Cogenus 3 case}
(Table 1.)
We fix the degree and the number of nodes to explain more clearly.
 The way to get a
general formula is exactly the same, but much more computation.
Let $d=5$ (degree 5 curves) and $\delta =3$ (3 nodes). There is
a technical reason to assume that $d > \delta$. It makes computation easy.
Then we have 11 points on $\pr^2$ and 6 points on
$\tilde{\pr^2}$. The formula of [V] gives $\text{N} (3,5)=7915$.
\noindent
\vspace{-3 mm}
\begin{enumerate}
 \item[(1)] Case A.  \par
   This case is easy. It's just $\text{N} (3,4)  =
   deg \ \overline{\VV (3,4)}$ = 675. Generally, the degree of this locus
   is $\text{N}_{3, d-1}$
 \item[(2)] Case B. \par
   $C_1 \cap \text{E} = C_2 \cap \text{E} = [4]$. This means that $C_1$ and E
   meet in 4 distinct points. Deg $\ C_1 = 4$, but the condition of $C_1$
   has 11 points + 2 nodes = 13 conditions, so one condition is needed.
   One of the four points on $C_1 \cap \text{E}$ gives this condition.
   So the number of all curves through 12 points (given 11 points + one
   axis point) and having two nodes is $\text{N} {(2,4)}$ (generally,
   $\text{N} {(2, d-1)}$). And
   the others on $\text{E} \cap C_1$ give three conditions on $C_2$. $C_2$
   has one node, so has one ruling. So the number of ways of choosing a ruling
   is 6 + 3 (generally $(d+1)+(d-2)$). So degree of Case B is
   $\text{N} (2,d-1) \times (2d-1)$.
 \item[(3)] Case C and Case D. \par
   The method is similar as that of B.
\end{enumerate}
\bigskip
 Table 1. (Calculation of N(3,5), and
          General Case N(3,$d$))
 \newline
 \noindent
 \begin{tabular}{|c|c|c|c|c|r|l|} \hline
 Case & \ M\  & N($C_1$) & N($C_2$) &
 \ \ $\pi$ on E & d=5 & General Case \\
 \hline
 A & 1 & 3 & 0 & [d-1] & 675 & $1 \times \text{N} {(3,d-1)} \times 1$  \\
 B & 1 & 2 & 1 & [d-1] & 2025 & $1 \times \text{N} {(2,d-1)} \times (2d-1)$  \\
 C & 1 & 1 & 2 & [d-1] & 756 & $1 \times \text{N} {(1,d-1)}
  \times (2d^2-5d+3)$\\
 D & 1 & 0 & 3 & [d-1] & 35 & $1 \times 1 \times (4d^3-24d^2+47d-30)/3$\\
 \hline
 E & 2 & 2 & 0 & [d-3,1] & 2020 & $2 \times \text{N} {([d-3,1],2,d-1)}
  \times 1$ \\
 \hline
 F & 2 & 1 & 1 & [d-3,1] & 1316 & $2 \times (6d^3-42d^2+90d-56)
 \times (2d-3)$ \\ \cline{6-7}
   &   &   &   &         & 200  & $2 \times (3d^2-12d+10) \times 2(d-3)$ \\
 \hline
 G & 2 & 0 & 2 & [d-3,1] & 24 & $2 \times 1 \times (4d^2-24d+32)$
 \\ \cline{6-7}
   &  &    &    &        & 60 & $2 \times 2(d-4) \times (2d^2-9d+10)$  \\
 \hline
 H & 3 & 1 & 0 & [d-4,0,1] & 405 & $3 \times \text{N} {([d-4,0,1],1,d-1)}
 \times 1$ \\
 \hline
 I & 3 & 0 & 1 & [d-4,0,1] & 54 & $3 \times 3(d-4) \times (2d-4)$
 \\ \cline{6-7}
   &   &   &   &           &  9 & $3 \times 1 \times 3(d-4)$      \\
 \hline
 J & 4 & 1 & 0 & [d-5,2] & 320 & $4 \times \text{N} {([d-5,2],1,d-1)}
 \times 1$ \\
 \hline
 K & 4 & 0 & 1 & [d-5,2] & 0 & $4 \times 2(d-4)(d-5)
 \times (2d-5)$ \\ \cline{6-7}
   &   &   &   &         & 0 & $4 \times 2(d-4) \times 2(d-5)$ \\
 \hline
 L & 4 & 0 & 0 & [d-5,0,0,1] & 16 & $4 \times 4(d-4) \times 1$\\
 \hline
 M & 6 & 0 & 0 & [d-6,1,1] & 0 & $6 \times 6(d-4)(d-5) \times 1$\\
 \hline
 N & 8 & 0 & 0 & [d-7,3] & 0 & $8 \times 4(d-4)(d-5)(d-6) \times 1$\\
 \hline
 Sum & & & & & 7915 & See *7 \\
 \hline
 \end{tabular}
 \newline
 *1. M means multiplicity.
 \newline
 *2. N($C_1$) means the number of nodes of $C_1$.
 \newline
 *3. N($C_2$) means the number of nodes of $C_2$.
 \newline
 *4. $\text{N} {(\delta,d)}$ = $deg \ \overline{\VV (\delta,d)}$.
 \newline
 *5. $\text{N} {([\pi],\delta,d)}$ = degree of the closure of
 the set of all curves $C$ of degree d having $\delta $
 \newline
 nodes
 and $C \cap \text{L}$ is {\it of type} $\pi$ to fixed line L.
 \newline
 *6. In General Case, M $\times$ A $\times$ B means that
  M is the multiplicity, A is the degree
 \newline
 of the locus of
 curves $C_1$ in $\pr^2$ and B is the degree of the locus of curves
 $C_2$ in $\tilde{\pr^2}$.
 \newline
 *7. N(3,d) = $9/2d^6-27d^5+9/2d^4+423/2d^3-229d^2-829/2d+525$.
 \par
 \newpage
\begin{enumerate}
 \item[(4)] Case E. \par
   Since $C_2$ is nonsingular, each point on $C_1 \cap \text{E} =
   \pi = [2,1]$ (generally, [$d-3,1$]) gives a condition on $C_2$,
    therefore the degree of
   curves $C_2$ is just one. $C_1$ is the curve with two nodes
    and one tangency
   to E. So $C_1$ has exactly 14 conditions (11 points + 2 nodes +
   1 tangency).
    Therefore, the degree of Case E is the same as the degree
   of the locus of all curves with two nodes and one tangency condition
   to fixed line. The method of calculating this degree
   is similar to that for the Severi variety $\VV (2,4)$, except when
   we degenerate $\pr^2$ to $\pr^2 \cup \tilde{\pr^2}$,
   we make the tangency on $C_1$ go to $\tilde{\pr^2}$, i.e,
 $$
 \begin{array}{ccccc}
 \pr^2 & \sim\rightarrow & \pr^2 & \cup & \tilde{\pr^2} \\
 \cup  &                 & \cup  &      & \cup   \\
 C_1   & \sim\rightarrow & C_{11} & \cup & C_{12} \\
 \cup  &                 & \cup  &      & \cup  \\
 11 \ \text{pts} & \sim\rightarrow & 7 \ \text{pts} & & 4 \text{pts} + 1
\  \text{tangency to} \ \text{E} \ .
 \end{array}
 $$
 Table 2 is the table of Case E{}.
 The interesting subcategory is Case H1. This case doesn't happen
 during calculation of the degree of the Severi variety $\VV (2,d-1)$.
  This case occurs when the
 curve $C_1$  degenerates to $C_{12}$ which contains a double ruling.
\end{enumerate}
\vspace{-4 mm}
 Table 2. (Calculation of $\text{N}{([d-3,1],2,d-1)}$)
 \newline
 \noindent
 \begin {tabular}{|c|c|c|c|c|r|l|} \hline
 Case & \ M\ & N($C_1$) & N($C_2$) &
  $\pi$ on E & d=5 & General Case \\
 \hline
 A1 & 1 & 2 & 0 & [d-2] & 126 & $1 \times \text{N} {(2,d-2)}
 \times 1 \times 2(d-2)$ \\
 B1 & 1 & 1 & 1 & [d-2] & 280 & $1 \times \text{N} {(1,d-2)}
 \times (2d-4) \times 2(d-3)$ \\
 C1 & 1 & 0 & 2 & [d-2] & 20 & $1 \times 1 \times (2d^2-11d+15)
 \times 2(d-4)$ \\
 \hline
 D1 & 2 & 1 & 0 & [d-4,1] & 432 & $2 \times \text{N} {([d-4,1],1,d-2)}
 \times 1 \times 2(d-2)$ \\
 \hline
 E1 & 2 & 0 & 1 & [d-4,1] & 64 & $2 \times 2(d-4)
 \times (2d-6) \times 2(d-3)$ \\
 \cline{6-7}
    &   &   &   &         & 16 & $2 \times 2(d-4) \times 2(d-3)$ \\
 \hline
 F1 & 3 & 0 & 0 & [d-5,0,1] & 54 & $3 \times 3(d-4) \times 1 \times 2(d-2)$ \\
 \hline
 G1 & 4 & 0 & 0 & [d-6,2] & 0 & $4 \times 2(d-4)(d-5) \times 1
 \times 2(d-2)$ \\
 \hline
 H1 & 1 & 0 & 2 & [d-4,2] & 10 & $2(d-1)+2(d-4)$ \\
 \hline
 Sum & & & & & 1010 & $(9d^5-90d^4+300d^3-327d^2-76d+190)$ \\
 \hline
 \end{tabular}
\begin{enumerate}
 \item[(5)] Case F. \par
 Since $C_1$ has 13 conditions (11 points + 1 node + 1 tangent condition),
 one more condition is needed. The point on $\text{E} \cap C_1$
 gives one more condition.
 Since the divisor $\text{E} \cap C_1$ has one multiple point, we have to
 divide two cases to fix one point of that divisor.\par
  \begin{enumerate}
   \item We fix an ordinary point. \par
    In this case, the condition of $C_1$ has 14 conditions. As in Case E, we do
    the whole thing again, we get the formula
     $6(d-3)^3+8(d-3)^2+2(d-3)+2(d-4)(2d-5) =
    (6d^3-42d^2+90d-56)$. The number of ways of choosing a ruling is
    $(2d-3)$, for the $(d+1)$ points in
    $\tilde{\pr^2}$ plus the $(d-4)$ points on axis $C_1 \cap E$.
    (We can't choose the fixed point on E
    which gives the condition on $C_1$ and tangent point.)
    So in this case, the degree
    is $(6d^3-42d^2+90d-56) \times (2d-3)$.
   \item We fix the tangent point. \par
    In this case, $C_1$ has 14 conditions.
     (11 points + 1 node + 1 tangent condition
    at fixed point.) As in Case E, we do the whole thing again, we get
    $3(d-3)^2+4(d-3)+(2d-5) = (3d^2-12+10)$.
    Since the tangent point is fixed on $C_1$,
    this point is not fixed on $C_2$. This gives the
    degree of the locus of rational smooth
    curves of $C_2$. This is two since this
    is the discriminant of the partition [0,1].
    And the number of ways of choosing a ruling is just $(d-3)$,
    for the $(d-3)$ points on
    $C_1 \cap$ E (as in previous subcase,
    we can't choose a tangent point. Also, we can't choose given inside points
    because there are no ordinary fixed points on $C_1 \cap E$).
    So in this case,
    the degree is $(3d^2-12d+10) \times 2(d-3)$.
  \end{enumerate}
 \item[(6)] Case G. \par
  Since $C_1$ has 12 conditions (11 points + 1 tangency on E),
  two conditions are needed. Two points on
  $C_1 \cap E$ give these conditions. But,
  as in Case F, the tangency condition gives us two ways to
  fix two points.
   \begin{enumerate}
    \item We fix one ordinary point and one tangent point. \par
     Since $C_1$ has full conditions, such a curve exists and is unique.
     The counting of degree of the locus of curves $C_2$ is a little tricky. If
     we choose both of the two rulings in $(d-4)$ points on
     the axis (as in Case F,
     we can't choose the fixed ordinary point and tangent point.), then
     the degree of the locus of rational smooth curves of $C_2$ is 4
     since this is the discriminant of
     the partition [1,1], and if we choose one ruling in axis point and one
     ruling in $(d+1)$ points in $\tilde{\pr^2}$,
     then the degree of locus of rational
     smooth curves is two (discriminant of [0,1]).
     So in this case, the degree is
     $4 \times \binom{d-4}{2} + 2(d+1)(d-4)$.
    \item We fix two ordinary points. \par
     The degree of the locus of curves $C_1$ is
      $2(d-4)$ (discriminant of $[d-5,1]$),
     and the number of ways of choosing a ruling is
     $\{ \binom{d+1}{2} +(d-5)(d+1) + \binom{d-5}{2} \}$. So in this
     case, the degree is $2(d-4) \times (2d^2-9d+10)$.
   \end{enumerate}
 \item[(7)] Case I, J, and K. \par
  The method to calculate these cases is similar as those of
  Case E, F, and G.
 \item[(8)] Case L. \par
  $C_2$ has no node, so all points on $E \cap C_1$ give conditions
  on $C_2$, so just one such curve exists.
  The degree of the locus of curves $C_1$ is
  just the discrimant of the divisor of $C_1 \cap {\text{E}} $.
   By Corollary 2.3, this is $4(d-4)$.
 \item[(9)] Case M and N. \ \
  The method is similar  as that of Case L.
\end{enumerate}

\subsection{Cogenus 4}

The way to calculate in the case of cogenus 4 is exactly the
same  as that of cogenus 3.
There are just more cases, so the only thing to do
is to be careful not to miss any. Case F involves
calculating the degree of the variety of all curves
with three nodes and one tangency to fixed line.
This counting method is the same as the case E of cogenus 3.
Table 3 and 4 are the calculations of low degrees.
Below each Table,  we give the general formula of the degree of these cases.
Compared with Table 1, Table 4 has new cases (Cases O1, P1, and Q1).
 These Cases correspond
to the cases having double rulings in $\tilde{\pr^2}$.
\par
Table 3. (Calculation of N(4,5), and
                         N(4,6))
\newline
\noindent
\begin{tabular}{|c|c|c|c|c|r|r|} \hline
\ Case\  & Multiplicity &  Nodes of $C_1$ & Nodes of $C_2$ &
\ \ \ $\pi$ on E & \ $d=5$ \ & \ $d=6$ \ \\
\hline
A & 1 & 4 & 0 & [d-1] & 666 & 36975 \\
B & 1 & 3 & 1 & [d-1] & 6075 & 87065 \\
C & 1 & 2 & 2 & [d-1] & 6300 & 39690 \\
D & 1 & 1 & 3 & [d-1] & 945  & 4032 \\
E & 1 & 0 & 4 & [d-1] & 15   & 70 \\
\hline
F & 2 & 3 & 0 & [d-3,1] & 4728 & 99160 \\
G & 2 & 2 & 1 & [d-3,1] & 10440 & 90708 \\
H & 2 & 1 & 2 & [d-3,1] & 1920 & 12800 \\
I & 2 & 0 & 3 & [d-3,1] & 0    & 224 \\
\hline
J & 3 & 2 & 0 & [d-4,0,1] & 2520 & 19170 \\
K & 3 & 1 & 1 & [d-4,0,1] & 1350 & 7002 \\
L & 3 & 0 & 2 & [d-4,0,1] & 0    & 252 \\
\hline
M & 4 & 2 & 0 & [d-5,2] & 1696   & 30240 \\
N & 4 & 1 & 1 & [d-5,2] & 0      & 5824 \\
O & 4 & 0 & 2 & [d-5,2] & 0      & 0 \\
P & 4 & 1 & 0 & [d-5,0,0,1] & 320 & 1344 \\
Q & 4 & 0 & 1 & [d-5,0,0,1] & 0   & 128 \\
\hline
R & 5 & 0 & 0 & [d-6,0,0,0,1] & 0 & 25 \\
\hline
S & 6 & 1 & 0 & [d-6,1,1] & 0     & 0 \\
T & 6 & 0 & 1 & [d-6,1,1] & 0     & 0 \\
\hline
U & 8 & 1 & 0 & [d-7,3] & 0 & 0 \\
V & 8 & 0 & 1 & [d-7,3] & 0 & 0 \\
W & 8 & 0 & 0 & [d-7,1,0,1] & 0 & 0 \\
\hline
X & 9 & 0 & 0 & [d-7,0,2] & 0 & 0 \\
\hline
Y & 12 & 0 & 0 & [d-8,2,1] & 0 & 0 \\
\hline
Z & 16 & 0 & 0 & [d-9,4] & 0 & 0 \\
\hline
Sum & & & & & 36975 & 437517 \\
\hline
\end{tabular}
\par
The general foumula  for $\text{d} \geq 5$ is :
$$
\text{N}(4,d) = 27/8d^8-27d^7+1809/4d^5-642d^4-2529d^3+37881/8d^2+
18057/4d-8865
$$

Table 4. (Calculation of $\text{N} {([2,1],3,4)}$, and
                         $\text{N} {([3,1],3,5)}$)
\newline
\noindent
\begin{tabular}{|c|c|c|c|c|r|r|} \hline
\ Case \ & Multiplicity & Nodes of $C_1$ & Nodes of $C_2$ &
\ \ $\pi$ \ on E & \ $d=4$ \ & \ $d=5$ \ \\
\hline
A1 & 1 & 3 & 0 & [d-2] & 90 & 5400 \\
B1 & 1 & 2 & 1 & [d-2] & 504 & 10800 \\
C1 & 1 & 1 & 2 & [d-2] & 240 & 2268 \\
D1 & 1 & 0 & 3 & [d-2] & 0 & 40 \\
\hline
E1 & 2 & 2 & 0 & [d-4,1] & 360 & 16160 \\
F1 & 2 & 1 & 1 & [d-4,1] & 672 & 7968 \\
G1 & 2 & 0 & 2 & [d-4,1] & 0 & 240 \\
\hline
H1 & 3 & 1 & 0 & [d-5,0,1] & 378 & 3240 \\
I1 & 3 & 0 & 1 & [d-5,0,1] & 0 & 324 \\
\hline
J1 & 4 & 1 & 0 & [d-6,2] & 0 & 2560 \\
K1 & 4 & 0 & 1 & [d-6,2] & 0 & 128 \\
L1 & 4 & 0 & 0 & [d-6,0,0,1] & 0 & 0 \\
\hline
M1 & 6 & 0 & 0 & [d-7,1,1] & 0 & 0 \\
\hline
N1 & 8 & 0 & 0 & [d-8,3] & 0 & 0 \\
\hline
O1 & 1 & 1 & 2 & [d-3,1] & 96 & 344 \\
P1 & 1 & 0 & 3 & [d-3,1] & 24 & 60 \\
Q1 & 2 & 0 & 2 & [d-5,2] & 0 & 48 \\
\hline
Sum & & & & & 2364 & 49580 \\
\hline
\end{tabular}
\newline
Here, $\text{N} {(\pi,3,d-1)}$ is the degree of the locus
 of all plane curves $C$ of degree
 \newline
 $d-1$
 having 3 nodes and $C \cap \text{L}$
is {\it of type $\pi$} to fixed line L.

The general formula for $\text{d} \geq 4$ is :
$$
\begin{array}{l}
\text{N} {([d-3,1],3,d-1)}
=9d^7-126d^6+603d^5-891d^4-1118d^3+3223d^2+416d-2416 \\
\text{N} ([d-5,2],2,d-1)
=9d^6-135d^5+750d^4-1785d^3+1249d^2+1188d-1116 \\
\text{N} ([d-4,0,1],2,d-1)
=27/2d^5-297/2d^4+549d^3-1341/2d^2-357/2d+495 \\
\text{N} ([d-6,1,1],1,d-1)
=18d^4-234d^3+1038d^2-1734d+720 \\
\text{N} ([d-5,0,0,1],1,d-1)
=12d^3-96d^2+220d-120 \\
\text{N} ([d-7,3],1,d-1)
=4d^5-76d^4-540d^3-1732d^2+2304d-720.
\end{array}
$$
\section{Degree of Severi polynomial}
\begin{defn}
A Severi polynomial is a polynomial $\text{N} (\pi,\delta,d)$
which is the degree of the locus of curves $C$ in
$\pr^2$ with degree d having $\delta$ nodes and
$C \cap \text{L}$ is {\it of type $\pi$} to fixed line L.
\end{defn}
\begin{rem}
The fact that $\text{N} (\pi,\delta,d)$ is a polynomial in d  is
obvious because of the recursion formula of Ran [R2].
\end{rem}
\begin{pro}
Fix $\delta$ and $\pi = [\ell_1,...,\ell_n]$.  Then
$$
\text{\rm deg N}(\pi,\delta,d) = 2\delta + \sum_{i=2}^{n} \ell_i.
$$
\end{pro}
{\it Proof} \ \ \
We use induction on $d$. $\text{N} (d)=\text{N} (\pi,\delta,d)$
 is the number of curves in
$\pr^2$ with degree $d$, having $\delta$ nodes, type $\pi$
conditions to fixed line and through given $n=(\binom{d+2}{2} - 1
- \delta - (\sum_{i=1}^{n} (i-1)\ell_i))$ points.
As in the previous section, we degenerate $\pr^2$ into $\pr^2 \cup
\tilde{\pr^2}$, and put $n_1=(\binom{d+1}{2} -1 - \delta)$ points into $\pr^2$
and $n_2=((d+1)-(\sum_{i=1}^{n} (i-1)\ell_i))$ points into $\tilde{\pr^2}$, and
degenerate the tangency conditions into $\tilde{\pr^2}$, i.e,
$$
\begin{array}{ccccc}
\pr^2 & \sim\rightarrow & \pr^2 & \cup & \tilde{\pr^2} \\
\cup  &                 & \cup  &      & \cup \\
C     & \sim\rightarrow & C_1   & \cup & C_2 \\
\cup  &                 & \cup  &      & \cup \\
n \ \text{points} & \sim\rightarrow & n_1 \ \text{points}
&  & n_2 \ \text{points}  +  \text{tangency conditions to L} \ .
\end{array}
$$
Then $\text{N} (\pi,\delta,d)$ is the sum of the Severi polynomials of all
limit components. So, if we prove that the degree of each polynomial
is less than or equal to $2\delta + \sum_{i=2}^{n} \ell_i$, we are done.
Each polynomial of a limit component is the product of
the polynomial of the upper part of that component (the set of the
locus $C_1$) and the polynomial of the lower part (the set of the
locus $C_2$). First, look at the degree of the polynomial
of the lower part. This degree is the sum of the number of rulings,
$\delta_2$, and the degree of the polynomial of the locus
 of the smooth curves with a divisor
{\it of type $\pi$} to fixed line, that is $\sum_{i=2}^{n} \ell_i$.
Second, look at the degree of the polynomial of the upper part.
$C_1$ has $\delta_1$ nodes and $C_1 \cap \text{E}$ is
{\it of type $\pi' =  [\ell_1',...,\ell_m']$}
satisfying the equation
\begin{equation}
\delta_1 + \sum_{i=1}^{m} (i-1)\ell_i'
= \delta - \delta_2
\end{equation}
Since $C_1$ has degree $d-1$, by the induction, the polynomial of the
upper part has degree at most $2\delta_1 + \sum_{i=2}^{m} \ell_i'$. So the
degree of the polynomial of each component is less than or equal to
$\delta_2 + \sum_{i=2}^{n} \ell_i + 2\delta_1 + \sum_{i=2}^{m} \ell_i'$.
 From the equation (2), we get
\begin{equation}
\delta_2 + \sum_{i=2}^{n} \ell_i + 2\delta_1 + \sum_{i=2}^{m} \ell_i'
\leq
2\delta + \sum_{i=2}^{n} \ell_i.
\end{equation}
 The equality holds just for
the component which is the locus of curves such that
  $C_1$ has $\delta$ nodes and
$C_2$ is a smooth rational curve.
\begin{cor}
Let $a_{2\delta}^{\delta}\ ( a_{2(\delta -1)}^{(\delta -1)})$ be the
coefficient of the degree $2\delta (2(\delta-1), respectively)$ term of the
Severi polynomial N($\delta$,d) (N($\delta-1$,d), respectively).
This is the leading coefficient of each polynomial by
the Proposition 4.2. Then
$$
 a_{2\delta}^{\delta} \times 2\delta
   = a_{2(\delta -1)}^{\delta -1} \times 6, \ \
 hence \ \ a_{2\delta}^{\delta} = {3^{\delta}\over{\delta !}}.
$$
\end{cor}
{\it Proof} \ \  \
In the proof of Proposition 4.2, equality holds for just
one limit component. The Severi polynomial of this component
is N($\delta, d-1$). Let's find the limit components of degree  one less than
the degree of Severi polynomial. By equation (3), there
exist exactly two such components. One such component, say $A_1$, consists of
curves such
that $C_1$ has ($\delta -1$) nodes and $C_2$ has one ruling, and
another component, say $A_2$, consists of
curves such that $C_1$ has ($\delta -1$) nodes
, $C_1 \cap \text{E}$ is {\it of type} [$d-3,1$] and $C_2$ is a
smooth rational curve. In Table 1(Table 3),
$ A_1$ is the component of Case B (Case B, respectively)
and $A_2$ is the component of Case E (Case F).
The coefficient of the degree $(2\delta -1)$ term of the polynomial of the
component $A_1$ is $2a_{2(\delta -1)}^{\delta -1}$. The coefficient of
the degree $(2\delta -1)$ term of the polynomial of the component $A_2$ is
$4a_{2(\delta -1)}^{\delta -1}$. So the sum is
$6a_{2(\delta -1)}^{\delta -1}$.
Therefore we get the following recursion
formula
$$
\text{N} (\delta,d) =
\text{N} (\delta,d-1) + 6a_{2(\delta -1)}^{\delta -1}d^{2\delta -1} \\
                    + \text{lower degree terms}.
$$
Integrating, we get $a_{2\delta}^{\delta} \times (2\delta)
                    = 6 \times a_{2(\delta -1)}^{\delta -1}$.

\begin{rem}
$\text{N} (7,d) = {243\over{560}} d^{14} + \text{\rm lower degree terms.}
$
\end{rem}
\begin{cor}
Let $a_{2\delta -1}^{\delta}$ be the coefficient of the degree
$(2\delta -1)$ term of the polynomial N($\delta$,d). Then
$$
a_{2\delta -1}^{\delta} = - (2\delta) \times a_{2\delta}^{\delta},
\ \ hence \ \ a_{2\delta -1}^{\delta} =
 {-2\times 3^{\delta}\over {(\delta -1)!}}.
$$
\end{cor}
{\it Proof} \ \ \
The way to prove this corollary is the same as that of corollary 4.3.
Let's find the components with the polynomial of
degree two less than the degree of Severi polynomial.
By equation (3), there exist exactly three such components.
First component, $B_1$, consists of curves such
 that $C_1$ has ($\delta -2$) nodes and $C_2$
has two rulings. Second component, $B_2$, consists of curves such
 that $C_1$ has ($\delta -2$)
nodes, $C_1 \cap \text{E}$ is {\it of type} [$d-3,1$] and $C_2$ has
one ruling. Third one, $B_3$, consists of curves such
 that $C_1$ has ($\delta -2$) nodes, $C_1 \cap \text{E}$ is
{\it of type} [$d-5,2$] and $C_2$ is smooth.
In Table 1(Table 3), $B_1$ is the component of Case C (Case C, respectively)
and $B_2$ is the component of Case F (Case G, respectively) and $B_3$ is the
component of Case J (Case M, respectively).
Calculation of the coefficients of the degree $(2\delta -2)$ term of those
polynomials give $2a_{2(\delta -2)}^{\delta -2},
8a_{2(\delta -2)}^{\delta -2}, 8a_{2(\delta -2)}^{\delta -2}$
respectively. So the sum is $18a_{2(\delta -2)}^{\delta -2}$. The polynomial
of the limit components $A_1$ and $A_2$ also
 have a term of this degree.
The coefficient of this degree term of polynomial of $A_1$ is $
-a_{2(\delta -1)}^{\delta -1} + 2(-2(\delta -1)
a_{2(\delta -1)}^{\delta -1} + a_{2\delta -3}^{\delta -1})$.
The coefficient of this degree term of polynomial of $A_2$ is
$-8a_{2(\delta -1)}^{\delta -1} +4(-4(\delta -1)
a_{2(\delta-1)}^{\delta -1} + a_{2\delta -3}^{\delta -1}) +
8a_{2(\delta-2)}^{\delta -2}
+16a_{2(\delta -2)}^{\delta -2}$.
Using $a_{2\delta -3}^{\delta -1} = -2(\delta -1)
a_{2\delta -2}^{\delta -1}$ (by induction on $\delta$)
and $a_{2(\delta -1)}^{\delta -1} \times 2(\delta -1)
= 6a_{2(\delta -2)}^{\delta -2}$
(by Corollary 4.3), we get that the coefficient of this degree term is
$-9a_{2(\delta-1)}^{\delta -1}(2\delta -1)$.
So, we get the following formula,
$$
\begin{array}{ll}
\text{N} (\delta,d) =&\text{N} (\delta,d-1)
+6a_{2(\delta -1)}^{\delta}d^{2\delta -1}
-9a_{2(\delta -1)}^{\delta -1}(\delta -2)d^{2\delta -2} \\
&+ \ \text{lower degree terms}
\end{array}
$$
Integrating, we get
$$
\text{N} (\delta,d) = a_{2\delta}^{\delta}d^{2\delta}
-(2\delta)a_{2\delta}^{\delta}d^{2\delta -1} + \text{lower degree terms.}
$$
\begin{rem}
$\text{N} (7,d) \ = \
 {243\over{560}}d^{14} - {243\over{40}}d^{13} + \text{\rm lower degree terms.}
$
\end{rem}

\noindent
Department of Mathematics, University of Durham,
Durham, UK

\noindent
Youngook.Choi@@durham.ac.uk

\end{document}